\documentclass[11pt,a4paper,leqno]{amsart}
\usepackage{amssymb}
\numberwithin{equation}{section}

\newtheoremstyle{theorem}{3pt}{3pt}%
{\it}
{}
{\bfseries}
{:}
{.5em}
{}
\theoremstyle{theorem}
\newtheorem{theorem}{Theorem}[section]

\newtheorem{proposition}[theorem]{Proposition}
\newtheorem{corollary}[theorem]{Corollary}
\newtheorem{lemma}[theorem]{Lemma}
\newtheorem{definition}[theorem]{Definition}

\newtheoremstyle{example}{3pt}{3pt}%
{}
{}
{\sc}
{:}
{.5em}
{}
\theoremstyle{example}

\newtheoremstyle{remark}{3pt}{3pt}%
{}
{}
{\sc}
{:}
{.5em}
{}
\theoremstyle{remark}
\newtheorem{remark}{Remark}[section]

\numberwithin{equation}{section}

\newcommand{\thismonth}{\ifcase\month\or
  January\or February\or March\or April\or May\or June\or
  July\or August\or September\or October\or November\or December\fi
  \space\number\year}
\makeatletter
\newcommand{\low}{\@ifnextchar^{}{^{\vphantom x}}}
\newcommand{\high}{\@ifnextchar_{}{_{\vphantom I}}}
\makeatother
\DeclareSymbolFont{script}{U}{eus}{m}{n}
\DeclareSymbolFontAlphabet{\mathscr}{script}
\DeclareMathSymbol{\EuWedge}{0}{script}{"5E}
\DeclareMathAlphabet{\mathrmsl}{OT1}{cmr}{m}{sl}
\newcommand{\rssymb}[2]{\newcommand{#1}{{\mathrmsl{#2}}}}
\newcommand{\calsymb}[2]{\newcommand{#1}{{\mathcal{#2}}}}
\newcommand{\bbsymb}[2]{\newcommand{#1}{{\mathbb{#2}}}}
\newcommand{\lieoper}[2]{\newcommand{#1}{\mathop
  {\mathfrak{#2}\null}\nolimits}}
\newcommand{\oper}[3][n]{\newcommand{#2}{\mathop
  {\mathrm{#3}\null}\ifx n#1\nolimits\else\limits\fi}}
\newcommand{\rsoper}[3][n]{\newcommand{#2}{\mathop
  {\mathrmsl{#3}\null}\ifx n#1\nolimits\else\limits\fi}}
\bbsymb\C{C} \bbsymb\F{F} \bbsymb\HQ{H}\bbsymb\N{N} \bbsymb\Q{Q}
\bbsymb\R{R} \bbsymb\U{U} \bbsymb\V{V} \bbsymb\W{W} \bbsymb\Z{Z}
\bbsymb\bbf{F} \bbsymb\bbk{K} \bbsymb\bbi{I} \bbsymb\bbl{L} \bbsymb\bbo{O}
\bbsymb\bbj{J}
\bbsymb\bby{Y}
\bbsymb\bbp{P}
\bbsymb\bba{A}
\calsymb\cA{A} \calsymb\cB{B} \calsymb\cC{C} \calsymb\cD{D} \calsymb\cE{E}
\calsymb\cF{F} \calsymb\cG{G} \calsymb\cH{H} \calsymb\cI{I} \calsymb\cJ{J}
\calsymb\cK{K} \calsymb\cL{L} \calsymb\cM{M} \calsymb\cN{N} \calsymb\cO{O}
\calsymb\cP{P} \calsymb\cQ{Q} \calsymb\cR{R} \calsymb\cS{S} \calsymb\cT{T}
\calsymb\cU{U} \calsymb\cV{V} \calsymb\cW{W} \calsymb\cX{X} \calsymb\cY{Y}
\calsymb\cZ{Z}
\renewcommand{\geq}{\geqslant} 
\oper\End{End} \oper\Hom{Hom}                    
\oper\Sym{Sym} \oper\Skew{Skew}
\oper\Aut{Aut}                                   
\oper\GL{GL} \oper\SL{SL}\oper\Symp{Sp}
\oper\CO{CO} \oper\On{O} \oper\SO{SO} \oper\Pin{Pin} \oper\Spin{Spin}
\oper\CU{CU} \oper\Un{U} \oper\SU{SU} \oper\PSU{PSU}
\rsoper\Diff{Diff} \rsoper\SDiff{SDiff}
\lieoper\der{der}                                
\lieoper\gl{gl} \lieoper\sgl{sl}\lieoper\symp{sp}
\lieoper\co{co} \lieoper\so{so} \lieoper\spin{spin}
\lieoper\cu{cu} \lieoper\un{u}  \lieoper\su{su}
\rsoper\Vect{Vect} \rsoper\Ham{Ham}
\def\la#1{\hbox to #1pc{\leftarrowfill}}
\def\ra#1{\hbox to #1pc{\rightarrowfill}}

\newcommand{\norm}[2][]{|\mkern-2mu|#2|\mkern-2mu|
  _{\lower1pt\hbox{${}_{#1}$}}}
\newcommand{\Norm}[2][]{\bigl|\mkern-3mu\bigr|#2\bigr|\mkern-3mu\bigr|
  _{\lower1pt\hbox{${}_{#1}$}}}

\rsoper\dimn{dim}                           
\rsoper\grad{grad}                          
\rsoper\kernel{ker}\rsoper\image{im}        
\rsoper\alt{alt}   \rsoper\sym{sym}         
\rsoper\Ad{Ad}     \rsoper\ad{ad}           
\rsoper\CoAd{CoAd} \rsoper\coad{coad}       
\rsoper\trace{tr}  \rsoper\trfree{tf}       
\rsoper\detm{det}                           
\rsoper\Vol{Vol}                            
\rsoper\divg{div}                           
\rsoper\sign{sign}                          
\rssymb\iden{id}                            
\rssymb\vol{vol}                            
\oper\Imag{Im}\oper\Real{Re}                
\newcommand{\sd}{{\raise1pt\hbox{$\scriptscriptstyle +$}}}
\newcommand{\asd}{{\raise1pt\hbox{$\scriptscriptstyle -$}}}
\newcommand{\sdasd}{{\raise1pt\hbox{$\scriptscriptstyle\pm$}}}
\newcommand{\asdsd}{{\raise1pt\hbox{$\scriptscriptstyle\mp$}}}

\rsoper\scal{scal}
\def\kahl/{k\"ahler}
\def\Kahl/{K{\"a}hler}
\begin{document}
\title[Lorentzian Sasaki-Einstein Metrics on Connected Sums of $S^{2}\times S^{3}$]
{Lorentzian Sasaki-Einstein Metrics on Connected Sums of $S^{2}\times S^{3}$}
\author[R.Gomez]{Ralph R. Gomez}
\address{Ralph R. Gomez: Department of Mathematics and Statistics,
Swarthmore College, Swarthmore, PA 19081, USA.}
\email{rgomez1@swarthmore.edu}

\date{\thismonth}
\begin{abstract}
Negative Sasakian manifolds can be viewed as Seifert-$S^1$ bundles
where the base orbifold $\mathcal{Z}^{orb}$ has $K_{\mathcal{Z}^{orb}}$ ample such that $c_{1}(D)$ is a torsion class, where $\mathcal{D}$ is the contact subbundle.
We use this framework to settle completely an open problem formulated by Boyer and Galicki in \cite{BG5} which asks whether or not
$\#kS^{2}\times S^{3}$ are negative Sasakian manifolds for all $k$. As a consequence of the affirmative answer to this problem,
there exists so-called Sasaki $\eta$-Einstein and Lorentzian Sasaki-Einstein metrics on these five-manifolds for all $k$ and moreover all of these can be realized as links of isolated hypersurface singularities defined by weighted homogenous polynomials. The key step is to construct infinitely many hypersurfaces in weighted projective space that contain
branch divisors $\Delta=\sum_{i}(1-\frac{1}{m_{i}})D_i$ such that the $D_i$ are rational curves.
\end{abstract}
\maketitle
\vspace{-2mm}

\bigskip

\section{Introduction}
K\"ahler geometry has an odd-dimensional counterpart known as Sasakian geometry. This branch of Riemannian geometry has
received a lot of attention within the last decade particularly because of its ability to construct many
different types of odd-dimensional manifolds that admit Einstein metrics of positive Ricci curvature.
Such manifolds have also come to play a major role in theoretical physics, most notably in testing the veracity of the
AdS/CFT conjecture (see \cite{Sparks} for recent survey).
Recall that a Sasakian manifold possesses a nowhere vanishing vector field $\xi$, called the Reeb vector field, which foliates the manifold in such a way that transversely the geometry is K\"ahler. This 1-dimensional foliation, denoted by $\mathcal{F}_{\xi}$, is sometimes called the characteristic foliation. By defining a kind of `transverse first Chern class' called the
\emph{basic first Chern class} denoted by, $c^{B}_{1}(\mathcal{F}_{\xi})$, one can
then try to explore Sasakian structures on manifolds for which this basic first Chern class vanishes or is represented by a positive or negative
$(1,1)$ form. For example, Sasaki-Einstein manifolds fall under the category of positive Sasakian geometry, that is, the case in which
$c^{B}_{1}(\mathcal{F}_{\xi})>0$. In this realm, there are already many such examples and in dimension five a lot more can be said in the direction of classification of Sasaki-Einstein manifolds. (See \cite{BG5} for detailed survey.)

\indent On the other hand, there is not much known about precisely which manifolds admit a negative Sasakian geometry ($c^{B}_{1}(\mathcal{F}_{\xi})<0$) other than numerous examples. A primary reason for this is the following: a negative Sasakian manifold can be viewed as a Seifert-$S^1$ bundle over a compact K\"ahler orbifold such that the orbifold canonical divisor is ample and so a complete understanding of negative Sasakian manifolds involves a much better understanding of orbifolds with ample orbifold canonical divisor and there is not a lot known here in the direction of classification. But, in light of a theorem due to J. Koll\'ar in dimension five \cite{Kol1}, a better understanding of negative Sasakian manifolds can be achieved. It is the purpose of this article to show that
$\#kS^{2}\times S^{3}$ are negative Sasakian manifolds for all $k$ ($k=0$ which corresponds to $S^{5}$, $k=1,2$ are already known to hold \cite{BGM}). This settles completely, in the affirmative, an open problem
formulated by C.Boyer and K.Galicki in \cite{BG5} where they asked whether or not all connected sums of $S^{2}\times S^{3}$ are negative Sasakian manifolds. Moreover, we have some consequences of this result. Negative Sasaki $\eta$-Einstein and Lorentzian Sasaki-Einstein metrics also exist on these manifolds.
A Sasaki $\eta$-Einstein metric is a Sasakian metric (see section 1) which satisfies a more general Einstein condition namely $Ric_g=\lambda g + \nu \eta \otimes \eta$ where $\eta$ is a contact 1-form, and $\lambda$, $\nu$ are constants.  Sasaki $\eta$-Einstein manifolds are transversely K\"ahler- Einstein and Sasaki $\eta$-Einstein metrics have constant scalar curvature equal to $2n(\lambda+1).$ These manifolds can be regarded as odd-dimensional analogs of K\"ahler-Einstein manifolds. Lorentzian Sasaki-Einstein metrics are pseudo-Riemannian Sasaki metrics of signature $(1,2n)$ that satisfy the Einstein condition $Ric_g=\lambda g.$
\\
\indent The main theorem of this article is:
\begin{theorem} Negative Sasakian structures exist on $\#kS^{2}\times S^{3}$ for all $k$. In particular, negative Sasaki $\eta$-Einstein
and Lorentzian Sasaki-Einstein metrics exist on these $5$-manifolds for all $k$.
\end{theorem}

Now, Sasaki-Einstein metrics were initially shown to exist on $\#k S^{2}\times S^{3}$ for $k=1,2$ in \cite{BG2}.
Later, the result was extended to
$k=1,\ldots,9$ (c.f.\cite{BG4},\cite{BGN1},\cite{BGN3}). J. Koll\'ar then showed in \cite{Kol3}, using different techniques, that Sasaki-Einstein metrics exist for all $k\geq 6$. Theorem $1.1$ may be viewed as a Lorentzian version of these results obtained by the authors. A notable difference is that our proof also exhibits all connected sums of $S^{2}\times S^{3}$ as links of isolated hypersurface singularities defined by weighted homogenous polynomials. (Compare also with \cite{BGN4}).\\
\indent The proof of Theorem $1.1$ implements a method developed in \cite{BG2} using links of isolated hypersurface singularities
defined by weighted homogenous polynomials. A key ingredient involves constructing infinitely many hypersurfaces in weighted projective space
with branch divisors such that the components of the branch divisors are rational curves. This will allow us to control the topology of the link by Theorem $50$ in \cite {Kol1}. \\
\indent The remaining three sections are organized as follows. The following two sections review some of the basic material of Sasakian geometry and then we
discuss the Lorentzian version. In the final section, we prove the main theorem which is a main result obtained in the author's Ph.D thesis \cite{Gmz}. We conclude by briefly discussing some important corollaries.

\section{Sasakian Geometry and Links of Isolated Hypersurface Singularities}
In this section, we review the basics of Sasakian geometry on links of isolated hypersurface singularities defined by
weighted homogeneous polynomials. See \cite{BG5} for details.

\begin{definition}
A Sasakian manifold $(M,g)$ is an odd-dimensional compact Riemannian manifold if any, hence all, of the
following equivalent conditions hold:

\indent i.) There exists a Killing vector field $\xi$ of unit length on $M$ so that the
tensor field $\Phi$ of type $(1,1)$, defined by $\Phi(X)= -\nabla_{X}\xi$ satisfies
the condition $$(\nabla_{X}\Phi)(Y)=g(X,Y)\xi-g(\xi,Y)X$$ for any pair of vector fields
$X$ and $Y$ on $M$.

\indent ii.) There exists a Killing vector field $\xi$ of unit length on $M$ so that
the Riemannian curvature satisfies the condition $$R(X,\xi)Y=g(\xi,Y)X-g(X,Y)\xi$$
for any pair of vector fields $X,Y$ on $M$.

\indent iii.) The metric cone
$(C(M),\overline{g})=(M \times\mathbb{R}_{+},dr^{2}+r^{2}g,d(r^{2}\eta))$
is K\"ahler, where $r$ is a coordinate on $\mathbb{R}_{+}$ and the complex structure $I$
is defined as $IY=\Phi Y+\eta(Y)\Psi$, where $\Psi=r\partial_r$ is the Euler vector field
such that $I\Psi=-\xi$ and $\eta$ is a contact 1-form defined below.
\end{definition}

The Reeb vector field gives rise to a $1$-dimensional foliation on the compact manifold $M^{2n+1}$ called the characteristic foliation $\mathcal{F}_{\xi}$.
The foliation is \textit{quasiregular} if for every $p\in M$ there exists a cubical neighborhood $U$ such that every leaf of the characteristic foliation intersects a transversal through $p$ at most a finite number of times $N(p)$.
(When $N(p)=1$, the foliation is said to be \textit{regular}.) Equivalently,
\begin{definition}
A Sasakian manifold is \emph{quasiregular} if all the leaves of the characteristic foliation
are compact, that is the leaves are circles.
\end{definition}

The $1$-form $\eta$ defined by $g(\xi,Y)=\eta(Y)$ is
a contact $1$-form dual to the Reeb vector field $\xi$, and satisfies the contact condition
$\eta\wedge (d\eta)^{n}\neq 0$ everywhere on $M.$ Define the contact subbundle $\mathcal{D}=
$ker$\eta$. The tangent bundle has an orthogonal splitting as$$TM={\mathcal{D}}
\oplus L_{\xi}$$ where $L_{\xi}$ is the trivial line bundle generated by $\xi$.
If $\Phi$ is restricted to $\mathcal{D}$, then it is an almost complex structure $J$ which
which turns out to be integrable and $(\mathcal{D},\Phi|_{\mathcal{D}}=J,d\eta,g|_{\mathcal{D}})$ is called a transverse K\"ahler structure. Furthermore, a Sasakian structure $\mathcal{S}$ on $M$ is characterized by the quadruple
of structure tensors $\mathcal{S}=(\xi,\eta,\Phi,g)$, where the Riemannian metric $g$ is naturally compatible
with the other structure tensors. (Sometimes it is common to write $(\mathcal{S},g)$ when one wishes to emphasize the metric and we will use this notation later.) A Sasakian manifold actually admits a whole family of Sasakian structures because
of so-called transverse homotheties. Specifically, one can deform a given Sasakian structure $\mathcal{S}=(\xi,\eta,\Phi,g)$ defined by:
\begin{equation} \xi'=a^{-1}\xi, \hspace{.4cm}  \eta'=a\eta, \hspace{.4cm} \Phi'=\Phi, \hspace{.4cm} g'=ag+(a^{2}-a)\eta \otimes \eta \end{equation}
where $a\in \mathbb{R}_+.$

\indent Now, the quasiregular Reeb vector field generates a locally free circle action on $M$. Then, we have an Seifert-$S^{1}$ bundle $S^{1}\rightarrow M\rightarrow \mathcal{Z}^{orb}$ where $\mathcal{Z}^{orb}$ is a complex K\"ahler orbifold with integral K\"ahler class $[\omega]\in H^{2}_{orb}(\mathcal{Z},\mathbb{Z})$ such that
the pull back of the K\"ahler form is $d\eta$. The process can be inverted in what is known as the Inversion Theorem \cite{BG1}. That is,
for a given integral K\"ahler class $[\omega]\in H^{2}_{orb}(\mathcal{Z},\mathbb{Z})$, one can form an $S^1$-Seifert bundle over $\mathcal{Z}^{orb}$
such that the total space $M$ will be a smooth Sasakian manifold where $\eta$ is a connection 1-form whose curvature is $\omega.$
(In general $M$ will be an orbifold but if all the uniformizing groups inject into $S^1$, then $M$ is a smooth manifold.)

The standard example of a quasiregular Sasaskian manifold is an odd dimensional sphere $S^{2n+1}$ with its standard metric. But, another rich source of examples of Sasakian
manifolds can be obtained by considering links of isolated hypersurface singularities
defined by weighted homogenous polynomials.
Consider the weighted $\mathbb{C}^{*}$ action on $\mathbb{C}^{n+1}$ defined by
$$\lambda : (z_{0},...,z_{n})\longmapsto (\lambda^{w_{0}}z_{0},...,\lambda^{w_{n}}z_{n})$$
where $w_{i}$ are the weights which are positive integers. We use the standard notation $\textbf{w}=(w_0,...,w_n)$ to denote a weight vector. In addition, assume that not all of the weights
have a common factor, that is to say $gcd(w_{0},...,w_{n})=1.$ Let us first recall a definition.
\begin{definition}
A polynomial $f\in \mathbb{C}[z_{0},...,z_{n}]$ is weighted homogenous if it
satisfies $$f(\lambda^{w_{0}}z_{0},...,\lambda^{w_{n}}z_{n})=\lambda^{d}f(z_{0},...,z_{n})$$
where the positive integer $d$ is the degree of $f$.
\end{definition}
Let the \emph{weighted affine cone}, $C_{f}$, be a hypersurface in $\mathbb{C}^{n+1}$
defined by the equation $f(z_{0},...,z_{n})=0$ such that the origin is an isolated singularity.
We can now define the link of the isolated hypersurface singularity, denoted by $L_{f}$, as $$L_{f}=C_{f}\cap S^{2n+1}$$ where
$C_{f}$ is the \emph{weighted affine cone} defined by the equation $f(z_{0},...,z_{n})=0$ and
$S^{2n+1}$ is viewed as the unit sphere in $\mathbb{C}^{n+1}$. By the Milnor fibration
theorem \cite{Milnor}, $L_{f}$ is a smooth $2n-1$ manifold which is $n-2$ connected.
It is well-known that a link of an isolated hypersurface singularity defined by a weighted homogenous polynomial admits a quasiregular Sasakian structure
\cite{BG4} inherited from the ``weighted Sasakian structure" on the sphere, denoted by $S_{\textbf{w}}^{2n+1}$. So, the
quotient of $S_{\textbf{w}}^{2n+1}$ by the weighted locally free circle action $S_{\textbf{w}}^{1}$ is the weighted projective space $\mathbb{P}(\textbf{w})$ and moreover
the quotient of $L_{f}$ by the $S_{\textbf{w}}^{1}$ action cuts out a hypersurface $\mathcal{Z}_{f}$ in $\mathbb{P}(\textbf{w})$ in weighted projective space.
This leads to the useful commutative diagram:
\begin{equation}
\begin{matrix}
L_{d} &&\longrightarrow&& S^{2n+1}_{\textbf{w}} \\
 \Big\downarrow \pi &&&&\Big\downarrow\\
\mathcal{Z}_{f} &&\longrightarrow&& \mathbb{P}(\textbf{w}),
\end{matrix}
\end{equation}
where the vertical arrows are Seifert-$S^1$ bundles and the horizontal rows are Sasakian and K\"ahler embeddings respectively \cite{BGN2}.
Since the circle action is locally free, this implies the weighted hypersurface is a complex orbifold.
In general, complex orbifolds can play a key role in quasiregular Sasakian geometry.
\begin{remark}In section 3, we will also use the notation $\mathcal{Z}_d$ or $L_d$ to alternatively emphasize the degree of the polynomial.
\end{remark}
When the codimension of the orbifold singular set $\Sigma^{orb}(\mathcal{Z})$ in $\mathcal{Z}$ is at least two, it is unnecessary to distinguish between the algebraic singularities on the
underlying algebraic variety and the orbifold singularities. However, we will find it more advantageous to consider codimension one singular sets
of the orbifold. In this situation, the orbifold singular set contains a divisor on the orbifold and such a divisor is said to be the branch divisor.
\begin{definition} A \emph{branch divisor} $\Delta$ of an orbifold is a $\mathbb{Q}-divisor $ on the underlying variety
of the form $$\displaystyle\sum_{i}(1-\frac{1}{m_{i}})D_{i} .$$ Each $D_{i}$ is a Weil divisor on
the underlying variety which lies in $\Sigma^{orb}(\mathcal{Z})$ and the \emph{ramification index} $m_{i}$ of $D_{i}$ is the gcd of the orders of the local uniformizing groups taken over all points of $D_{i}$.
\end{definition}

Using weighted projective spaces, it is straightforward to find the orbifold singular locus
$\Sigma^{orb}(\mathbb{P}(\textbf{w}))$. Determination of this locus will be crucial in section 3 since the orbifold singular locus will
be a branch divisor of the orbifold.
Let $P$ be any proper subset of $\{w_{0},...,w_{n}\}$ and
let $S_{P}$ be defined as follows:
$$S_{P}=\{[z_{0}:...:z_{n}]\in \mathbb{P}(\textbf{w})\mid z_{j}=0~\forall w_{j} \in P \}.$$
Let$P^{c}=\{w_{i_1},...,w_{i_k}\}$ be the compliment of $P$ in $\{w_{0},...,w_{n}\}$. Then $S_{P}$ is in
$\Sigma^{orb}(\mathbb{P}(\textbf{w}))$ iff \\
$gcd(w_{i_{1}},...,w_{i_{k}})>1$ where the local uniformizing
groups for points in $S_{P}$ are the cyclic groups of the form $\mathbb{Z}_{gcd(w_{i_{1}},...,w_{i_{k}})}$.
Furthermore, weighted projective spaces with codimension two or more orbifold singular set can be numerically described as
follows (\cite{BG2}, \cite{BG5},\cite{Flet}):
\begin{definition}
A weighted projective space $\mathbb{P}(w_{0},...,w_{n})$ is
well-formed if $$gcd(w_{0},\ldots,\widehat{w}_{i},\ldots,w_{n})=1$$ for all $i=0,...,n$ and the hat indicates deletion of the $i^{th}$
weight. A hypersurface $\mathcal{Z}_f$ in $\mathbb{P}(\textbf{\emph{w}})$ is well-formed if in addition it contains no codimension 2 singular
stratum of $\mathbb{P}(w_{0},...,w_{n})$.
\end{definition}
Therefore, a \emph{non well-formed} weighted projective space satisfies the property that there exists
an $i$ such that $gcd(w_{0},\ldots,\widehat{w}_{i},\ldots,w_{n})>1$. In such a situation, the orbifold singular set
becomes a branch divisor of the orbifold.
We will focus almost exclusively in the case in which a branch divisor exists on the underlying variety.

\begin{remark}
Sometimes the common notation $(\mathcal{Z},\Delta)$ shall be implemented instead of $\mathcal{Z}^{orb}$ to emphasize
the presence of a branch divisor.
\end{remark}

Before we close the section, we must make sure that our hypersurfaces in weighted projective space give rise to links with isolated hypersurface singularities. Such conditions are well known and these are the \textit{quasismoothness} conditions (e.g. see \cite{Flet}):
\begin{theorem}\emph{(Quasismoothness Condition)}
Let $\mathcal{Z}_{f}$ be a hypersurface in $\mathbb{P}(w_{0},...,w_{3})$ where $f$ has no linear terms. Then
$\mathcal{Z}_{f}$ is quasismooth if and only if all of the three following conditions holds:
\begin{enumerate}
\item For each $i=0,1,2,3$, there is a $j$ and a monomial $z_{i}^{m_{i}}z_{j}\in \mathcal{O}(d)$.
The case $j=i$ is permitted.
\item If $gcd(w_{i},w_{j})>1$ then there exists a monomial $z_{i}^{b_{i}}z_{j}^{b_{j}}\in \mathcal{O}(d)$.
\item For every $i$,$j$ either there is a monomial $z_{i}^{b_{i}}z^{b_{j}}_{j}\in \mathcal{O}(d)$ or
there are monomials $z_{i}^{c_{i}}z^{c_{j}}_{j}z_{k}$ and $z_{i}^{d_{i}}z^{d_{j}}_{j}z_{l}\in \mathcal{O}(d)$
with $\{k,l\}\neq\{i,j\}$.
\end{enumerate}
\end{theorem}

\section{Sasaki $\eta$-Einstein and Lorentzian Sasaki-Einstein Structures}

Our aim in this section is to review the definition of \textit{negative}, \textit{positive}, and \textit{null} Sasakian structures on manifolds.
The definition depends upon the \textit{basic first Chern class} so the key ideas of basic cohomology will be recalled.
Intuitively, one can think of basic cohomology as the de Rham cohomology of the leaf space of the foliation (c.f. \cite{BG5},\cite{Mol},\cite{Tond}).
Let us briefly recall some of the key ideas of basic cohomology. An important aspect of basic cohomology in the context of Sasakian geometry
is that even though $[d\eta]_{B}$ is a trivial class in $H^{2}(M,\mathbb{R})$, it is a nontrival class in the basic cohomology group $H^{2}_{B}(\mathcal{F}_{\xi})$.
\begin{definition}
A smooth p-form $\alpha$ is basic if $$\xi\rfloor\alpha=0, \hspace{1cm}  \mathcal{L}_{\xi} \alpha=0.$$      
\end{definition}
Note that the exterior derivative of a basic form is again basic. Just as in de Rham cohomology,
one can define a basic de Rham complex for basic forms and this gives
rise to the basic cohomology ring denoted by $H^{*}_B(\mathcal{F}_{\xi})$. For compact $M$,
one can obtain a long exact basic cohomology sequence which relates to the usual deRham cohomology \cite{Mol}:
\begin{equation}
\cdots\longrightarrow H^{p}_{B}(\mathcal{F}_{\xi})\stackrel{\iota_{*}}{\longrightarrow}H^{p}(M,\mathbb{R})\stackrel{j_{p}}\longrightarrow
H^{p-1}_{B}(\mathcal{F}_{\xi})\stackrel{\delta}\longrightarrow H^{p+1}_{B}(\mathcal{F}_{\xi})\longrightarrow\cdots
\end{equation}
where $\iota_{*}$ is the natural inclusion, and $\delta$ is the connecting homomorphism defined by $\delta[\alpha]=[d\eta\wedge \alpha]=[d\eta]\cup [\alpha]$ and $j_{p}$ is defined by composing the map induced by $\xi\rfloor$ and the isomorphism involved in $H^{r}(M,\mathbb{R})\approx
H^{r}(M,\mathbb{R})^{T}$. The torus $T$ is the closure of the flow of the Reeb vector field and $H^{r}(M,\mathbb{R})^{T}$
means the $T$-invariant cohomology defined from the $T$-invariant r-forms $\Omega^{r}(M)^{T}$.
Moreover, since a Sasakian manifold is transversely K\"ahler, there are also the transverse Dolbeault cohomology groups denoted by $H_{B}^{p,q}(\mathcal{F}_{\xi})$.\\
\indent Because the contact subbundle $\mathcal{D}$ is complex, the first Chern class $c_{1}(\mathcal{D)}$
may be considered. Even more, $c_{1}(\mathcal{D)}$ can be represented by a real basic form $\rho_{B}$
in $H_{B}^{1,1}(\mathcal{F}_{\xi})\subset H^{2}_{B}(\mathcal{F}_{\xi})$. This $(1,1)$ form, as a class in $H_{B}^{2}(\mathcal{F}_{\xi})$, is denoted by $c^{B}_{1}(\mathcal{F}_{\xi})$
and is called the \textit{basic first Chern class}. We are now ready for:
\begin{definition}
A Sasakian structure on $M$ is \emph{negative} \emph{(positive)} if the basic first Chern class
$c^{B}_{1}(\mathcal{F}_{\xi})$ is represented by a negative(positive) definite $(1,1)$-form. In either         
case, the Sasakian structure is said to be \emph{definite} and \emph{indefinite} otherwise. When
$c^{B}_{1}(\mathcal{F}_{\xi})=0$ the Sasakian structure is called \emph{null}.
\end{definition}
There is a natural relationship between the orbifold first Chern class
and the basic first Chern class for a quasiregular Sasakian manifold. Specifically, on the base K\"ahler orbifold one can define
a Ricci form and so the pull back of this is the transverse Ricci form. This yields
\begin{lemma}\emph{\cite{BG5}}
Let $M$ be a quasiregular Sasakian manifold $S^{1}\longrightarrow M\stackrel{f}\longrightarrow (\mathcal{Z},\Delta)$.
Then the following relation holds: $$c^{B}_{1}(\mathcal{F}_{\xi})=f^{*}c_{1}^{orb}(\mathcal{Z},\Delta).$$
\end{lemma}
Therefore, a positive Sasakian structure on a manifold $M$ can be viewed as a smooth Seifert-$S^{1}$ bundle
$f:M\longrightarrow (\mathcal{Z},\Delta)$ such that $-K_{\mathcal{Z}^{orb}}$ is ample. $M$ is then said to have a\emph{ transverse
Fano structure}. Similarly, a negative Sasakian
manifold $M$ can be thought of as a smooth Seifert-$S^1$ bundle $f:M\longrightarrow (\mathcal{Z},\Delta)$ such that $K_{\mathcal{Z}^{orb}}$ is ample. In
this case, $M$ is said to have a \emph{transverse Aubin-Calabi-Yau} structure. Finally, a null Sasakian structure on a manifold $M$
is sometimes called a \emph{transverse Calabi-Yau} structure. Together with (2.2) and Lemma 3.3, finding negative Sasakian structures on links
of isolated hypersurface singularies defined by weighted homogenous polynomials gets translated into finding hypersurfaces in weigthed projective space with
orbifold first chern class negative or, equivalently, ample orbifold canonical class.

In the case $M$ admits a positive Sasakian structure for the basic first Chern
class proportional to $[d\eta]_{B}$,
then $M$ also admits a Sasakian metric of positive Ricci curvature \cite{BGN4}, and under additional conditions an
Einstein metric. What about the case in which $M$ admits a negative Sasakian structure? Before we turn to the geometry of this situation, we need to recall the $\eta$-Einstein condition:
\begin{definition}
A Sasakian structure $\mathcal{S}=(\xi,\eta,\Phi,g)$ on $M^{2n+1}$, $n>1$, is Sasaki $\eta$-Einstein if
$$Ric_{g}=\lambda g + \nu \eta\otimes\eta$$ where the constants $\lambda,\nu$
satisfy $\lambda + \nu = 2n.$
\end{definition}
\begin{remark}
Sasaki $\eta$-Einstein manifolds are of constant scalar curvature $2n(\lambda +1).$

\end{remark}
Let us now exploit a bit further, the transverse geometry of a Sasakian manifold under this additional $\eta$-Einstein hypothesis.
The transverse Ricci tensor $Ric_{g}^{T}$ of $g_{T}$ is the Ricci tensor of the transverse metric denoted by
$g_{T}$ and relates to the full Ricci tensor $Ric_{g}$ by $$Ric(X,Y)=Ric_{T}(X,Y)-2g(X,Y)$$
for any smooth sections $X$,$Y$ in $\mathcal{D}$. Further, we have $$\rho_{g}(X,Y)=Ric_{g}(X,\Phi Y),
 \hspace{.5cm} \rho_{g}^{T}(X,Y)=Ric^{T}_{g}(X,\Phi Y)$$
where $\rho_{g}$ and $\rho_{g}^{T}$ is the Ricci form and the transverse Ricci form respectively.
Just as in the K\"ahler case, the basic first Chern class can be represented by the transverse Ricci form.
Using the $\eta$-Einstein condition and the above relations, gives:
$$c_{1}^{B}(\mathcal{F_{\xi}})=(\lambda + 2) [d\eta]_{B}.$$ Since $[d\eta]_B$ is a nontrival class, we see that compact Sasaki $\eta$-Einstein manifolds can be arranged according to $\lambda<-2$ (negative $\eta$-Einstein), $\lambda=-2$ (null $\eta$-Einstein), and
$\lambda>-2$ (positive $\eta$-Einstein) and it is because of this that Sasaki $\eta$-Einstein manifolds are regarded
as the odd-dimensional analogs of compact K\"ahler-Einstein manifolds.

Because $M$ is transversely K\"ahler, there is a `transverse Calabi-Yau' problem
also and by the work of El Kacimi-Alaoui in \cite{ElKacimi}, there are no obstructions to solving the `transverse Monge-Ampere' equations when
the basic first Chern class is negative or null and so we have transverse versions of the celebrated theorems of Yau \cite{Yau} and Aubin \cite{Aubin}.
This is rephrased as
\begin{theorem}\cite{BGM}
Let $c_1^{B}(\mathcal{F}_{\xi})$ be proportional to $[d\eta]_B$. If $c_1^{B}(\mathcal{F}_{\xi})$ is zero or represented by a negative definite $(1,1)$ form, then there exists a Sasaki $\eta$-Einstein
structure with $\eta$-Einstein metric $g$ on $M$ with $\lambda=-2$  and $\lambda <-2$ in the other case.
\end{theorem}
For the case of compact five-manifolds, considerably more can be said because of a well-known classification theorem of Smale \cite{Sm}. Therefore, by
specializing to links of isolated hypersurface singularities, the topological data becomes encoded in the weights and degree of the weighted homogenous polynomial. Now, links of isolated hypersurface singularities are compact, simply-connected smooth $2n-1$ dimensional manifolds by the Milnor fibration theorem (for our case $n=3$)\cite{Milnor}.
In a natural way, these manifolds admit Sasakian structures (see for example \cite{BG2}.) Moreover, by Theorem 11.8.4 of \cite{BG5}, $c_1(\mathcal{D})=0$ and this is equivalent to the notion that $c_{1}^{B}$ is proportional to $[d\eta]_{B}$. Since $w_{2}(M)$ is the mod $2$ reduction of $c_{1}(\mathcal{D})$, we have that $w_{2}(L_f)=0$ and hence $L_{f}$ is spin. By Smale's classification theorem, any compact, simply-connected, spin five-manifold
for which $H_{2}(M,\mathbb{Z})$ has no torsion is diffeomorphic to $\#k S^{2}\times S^{3}$ where $k$ is the second Betti number of $M$.
Lastly, using an orbifold adjunction formula on $(\mathcal{Z},\Delta)$, a link $L_f$ admits a negative Sasakian structure if $d-\sum_{i}w_i>0$ and this was essentially established in \cite{BGK} (see also \cite{BG5}). A point worth emphasizing is that in the non well-formed cases, the canonical divisor and the orbifold canonical divisor do not coincide and therefore establishing an orbifold adjunction formula requires one work with the orbifold canonical divisor $K_{\mathcal{Z}^{orb}}$. Summarizing, we have
\begin{theorem}
Let $X_f$ be a quasismooth hypersurface in the weighted projective space \\$\mathbb{P}(w_0,w_1,w_2,w_3)$ and let $L_f$
be the corresponding link over this hypersurface. Then $c_{1}^{B}$ is proportional to $[d\eta]_{B}$ and  $L_f$ is a compact, simply-connected, spin Sasakian manifold of dimension five. Moreover, the following hold:\\
1.) If $H_{2}(L_{f},\mathbb{Z})$ has no torsion, then $L_f$ is diffeomorphic to $\#k S^{2}\times S^{3}$ where
$k$ is $b_2(L_f)$.\\
2.) If $d-\sum_{i}w_i >0$, then $L_f$ admits a negative Sasaki $\eta$-Einstein structure.
\end{theorem}
Another metric which can arise from a negative Sasakian structure is a Lorentzian
Sasaki-Einstein metric on $M$. Recall that for $\lambda >-2$ Tanno \cite{Tanno} showed that given a quasi-regular Sasaki $\eta$-Einstein manifold,
by applying a transverse homothety one can obtain a new Sasakian metric which is Einstein. We shall use a Lorentzian version of this result to obtain a Lorentzian Sasaki-Einstein metric. But first we need
\begin{definition} Let $(M^{2n+1},g)$ be a Lorentzian manifold of signature $(1,2n)$ where the
nowhere vanishing vector field $\xi$ is a time-like Killing vector field such that
$g(\xi,\xi)=-1$. We say $M$ is Sasakian if $\Phi(X)=-\nabla_{X}\xi$ satisfies the condition $(\nabla_{X}\Phi)(Y)=g(X,Y)\xi +g(\xi,Y)X$ and Lorentzian
Sasaki-Einstein if $g$ is Einstein.
\end{definition}

Of course, we can alternatively define everything in terms of the cone and define $M$ to be
$\emph{Lorentzian Sasakian}$ if $C(M)=(\mathbb{R}_{+},-dt^{2}+r^{2}g,d(r^{2}\eta))$
is psuedo-K\"ahler and \emph{Lorentzian Sasaki-Einstein} with Einstein constant $-2n$ if the cone is pseudo-Calabi Yau.

Now we can state the key proposition which connects a negative Sasaki $\eta$-Einstein structure on $M$ with a
Lorentzian Sasaki-Einstein structure.
\begin{proposition}\emph{\cite{BGM}}
Let $(M,\xi,\eta,\Phi,g)$ be a compact quasiregular Sasaki $\eta$-Einstein manifold
of dimension $2n+1$ and let $\lambda <-2$. Then $M$ admits a Lorentzian Sasaki-Einstein
structure such that
$$\xi^{'}=\frac{1}{a}\xi, \hspace{1cm} -g^{'}=ag+a(a-1)\eta \otimes \eta$$
where $a=\frac{\lambda +2}{2+2n}$.
\end{proposition}
\begin{proof} The proof is a straightforward calculation based on the fact that the Ricci tensor under transverse homothety
is related by
\begin{equation}
Ric_{g'}=Ric_g-2(a-1)g+(a-1)(2n+2+2na)\eta \otimes \eta.
\end{equation}
\end{proof}
In particular, by Theorem 3.5, any negative Sasakian manifold such that $c^{B}_{1}(\mathcal{F}_{\xi})$ is a multiple of $[d\eta]_{B}$
admits an $\eta$-Einstein metric where $\lambda < -2$ and by the above proposition, there also exists a Lorentzian Sasaki-Einstein
structure on M, i.e.
\begin{corollary}
Every negative Sasakian manifold such that $c^{B}_{1}(\mathcal{F}_{\xi})$ is a multiple of $[d\eta]_{B}$ admits a Lorentzian Sasaki-Einstein structure.
\end{corollary}

\section{Proof of Theorem $1.1$}

Before we begin the proof of our main theorem, we state a theorem of Koll\'ar \cite{Kol1} in dimension five which is crucial to the proof of our main theorem.

\begin{theorem} Let $M^{5}$ be a compact, simply-connected quasiregular Sasakian 5-manifold with corresponding base orbifold
$(\mathcal{Z},\sum(1-\frac{1}{m_i})D_i)$.
Then $$H_{2}(M,\mathbb{Z})=\mathbb{Z}^{k}\oplus\sum(\mathbb{Z}/m_{i})^{2g(D_{i})}$$ where $g(D_i)$ denotes the genus of the branch divisor.
\end{theorem}
The idea of the proof of Theorem $1.1$ is the following:
We construct hypersurfaces in weighted projective space with branch divisors whose components are rational curves. The links over these will therefore have no torsion by Theorem 4.1. Then according to Theorem 3.6, we need to check that the second Betti number is arbitrary and we must verify that these Sasakian links are indeed negative.\\
\indent Because we need data about the genus of a weighted curve, the following proposition will be used throughout:
\begin{proposition}\cite{Flet}
Let $C_{d}$ in $\mathbb{P}(w_{0},w_{1},w_{2})$ be a nonsingular curve and $d$ the degree of $C.$ Then the genus $g$ is given by
\begin{equation}2g(C_d)=\frac{d^{2}}{w_{0}w_{1}w_{2}}-d\sum_{i<j}\frac{gcd(w_{i},w_{j})}{w_{i}w_{j}}+
\sum_{i}\frac{gcd(d,w_{i})}{w_{i}}-1.\end{equation}
\end{proposition}

Also, since we need to compute the second Betti numbers of links we will need to implement an algorithm \cite{MilnorOrlick} devised by Milnor and Orlik
for computing Betti numbers of links.

The method for determining the second Betti number of the link entails unraveling
the characteristic polynomial of the monodromy map, called the Alexander polynomial, which arises in the exact sequence
$$0\longrightarrow H_{n}(L_{f},\mathbb{Z})\longrightarrow H_{n}(F,\mathbb{Z})\stackrel{\mathbb{I}-h_{*}}\longrightarrow H_{n}(F,\mathbb{Z})
\longrightarrow H_{n-1}(L_{f},\mathbb{Z})\longrightarrow 0$$ where $F$ is a fiber occurring in the well-known
Milnor Fibration Theorem and $h_{*}$ is the monodromy map which arises from the action of a generator
of $\pi_{1}(S^{1})$ on the homology of the fiber. The idea is that the number of factors of
$t-1$ which occur in the Alexander polynomial $\Delta(t)=det(t\mathbb{I}-h_{*})$ is precisely the
the Betti number of the link. Milnor and Orlik \cite{MilnorOrlick} computed this in terms of the degree and weights of the
weighted homogenous polynomial by understanding an associated divisor of the Alexander polynomial.
For any monic polynomial $f$ with roots $\alpha_{1},...,\alpha_{k}$ in $\mathbb{C}^{*}$, construct
its divisor $$ divisor (f)=<\alpha_{1}>+...+<\alpha_{k}>$$ as an element in the integral ring $\mathbb{Z}[\mathbb{C}^{*}]$.
Then the divisor of the Alexander polynomial is$$divisor \Delta(t)=1+\sum a_{j}\Lambda_{j}$$ where $\Lambda_{j}=div(t^{j}-1).$
So now form the irreducible rational weights $\frac{d}{w_i}=\frac{u_i}{v_i}$. The second Betti number of the link is then given by:
\begin{equation}
b_{2}(L_{f})=\displaystyle \prod_{i}(\frac{\Lambda_{u_i}}{v_i}-1)=1+\sum_{j}a_{j}.
\end{equation}
where the $\Lambda$'s adhere to the relations $\Lambda_{a}\Lambda_{b}=gcd(a,b)\Lambda_{lcm(a,b)}.$

The proof of the main Theorem $1.1$ can be broken up into three cases. Case I corresponds to the case in which
$b_2\geq 8$ and is even. Case II corresponds to the case $b_2$ is odd and $\geq3$.  Case III covers the remaining
cases when $b_2=4,6$. The cases for $b_2=1,2$ were already established in \cite{BGM}. In this section, we use the notation
$L_d$ and $\mathcal{Z}_d$ to indicate the degree of the particular weighted homogenous polynomials. Now that we finally have all of the tools in place,
we may begin the proof.

\noindent \textit{Proof of Theorem 1.1:}\\
\textbf{Case I}\\
Consider the hypersurface $\mathcal{Z}_{4k}=(x^{4}+y^{2}+z^{k}+t^{k}=0)\subset \mathbb{P}(k,2k,4,4)$ of degree $d=4k$
such that $gcd(2,k)=1$. $\mathcal{Z}_{4k}$ is easily seen to be quasismooth and it is also non well-formed with
ramification index $m=2$. Since this hypersurface is non well-formed, there
exists a branch divisor on the hypersurface defined by $C=(x=0)$ of degree $2k.$
In particular $C=(y^{2}+z^{k}+t^{k}=0)\subset \mathbb{P}(2k,4,4)=\mathbb{P}(k,2,2)$ so now we apply
the weighted genus formula (4.1) which yields$$2g(C)=k-2k\left(\frac{1}{2k}+\frac{1}{2}+\frac{1}{2k}\right)+2=0.$$
Therefore, the orbifold $(\mathcal{Z}_{4k},\Delta=\frac{1}{2}C)$ has one curve in the branch divisor and this curve
is rational, since its genus is zero. The corresponding link $L_{4k}$ of this hypersurface
therefore is a negative Sasaki $\eta$-Einstein manifold as soon as $k\geq 9$ by Theorem 3.6 and by Theorem 4.1 $H_2(L_f,\mathbb{Z})$
has no torsion. Now we compute the second Betti number of the link:
\begin{align*}
\textrm{divisor}\Delta(t)=(k-2)\Lambda_{4k}-(k-2)\Lambda_{2k}+(k-2)\Lambda_{k}+\Lambda_{4}-\Lambda_{2}+1\\
\end{align*}
and hence $b_{2}(L_{4k})=k-1$. Finally, by Theorem 3.6 and Corollary 3.9, we conclude
$L_{4k}$ is diffeomorphic to $\#(k-1)S^{2}\times S^{3}$ and it also admits Lorentzian Sasaki-Einstein metrics,
where $k\geq9$, $k$ odd. This completes case I.

\noindent\textbf{Case II}\\
In this case the following hypersurface of degree $d=4(4k+1)(4k+3)$ is considered:
$$\mathcal{Z}_{d}=(x^{4}+y^{8k+2}+z^{4k+1}t+t^{2k+1}z=0)\subset\mathbb{P}((4k+1)(4k+3),2(4k+3),4(4k+1),8(4k+1))$$
where $k\geq1$. This hypersurface is of Type III in the list given by Yau-Yu
which is quasi-smooth \cite{YauYu} and so $L_{f}$ is the link of an isolated
singularity at the origin. As with case I, the aim is to show that this hypersurface
contains curves in its branch divisor which have genus equal to zero. Note that $\mathcal{Z}_{4(4k+1)(4k+3)}$
is non well-formed so there exists a branch divisor. However, this branch divisor
has two curves: $D_1$ defined by $x=0$ and $D_2$ defined by $y=0$ with corresponding ramification indices
$m_1=2$ and $m_2=4k+1$ respectively. The curve $D_1$ of degree $2(4k+1)(4k+3)$ is of the form $D_{1}=y^{2(4k+1)}+z^{4k+1}t+t^{2k+1}z$ so that
$$D_{1}\subset\mathbb{P}(2(4k+3),4(4k+1),8(4k+1))=\mathbb{P}(4k+3,2(4k+1),4(4k+1)).$$
Applying the weighted genus formula (4.1) to $D_1$ gives
\begin{align*}
2g(D_{1})=\frac{4k+3}{2}-1-\frac{4k+3}{2}-\frac{1}{2}+2-\frac{1}{2}=0.\\
\end{align*}
Therefore, $D_1$ is a rational curve. A nearly identical calculation shows that the curve
$$D_2=x^4+z^{4k+1}t+t^{2k+1}z\subset\mathbb{P}(4k+3,4,8)$$
of degree $4(4k+1)$ is a rational curve as well.
Therefore, the orbifold $$(\mathcal{Z}_{4(4k+1)(4k+3)},\Delta=\frac{1}{2}D_1+\frac{k}{k+\frac{1}{4}}D_2)$$
has only rational curves in its branch divisor. Thus the corresponding link $L_{4(4k+1)(4k+3)}$ admits a negative Sasaki $\eta$-Einstein structure
since $d-\sum_i w_{i}=8k(8k+2)-(4k+3)^{2}>0$ for $k\geq1$ and Theorem 3.6. Furthermore the corresponding link has no torsion in $H_{2}(L_{f},\mathbb{Z})$ by Theorem 4.1. Next we compute the second Betti number and we obtain:\begin{align*}
divisor\Delta(t)&=(\Lambda_{4}-1)(\Lambda_{2(4k+1)}-1)(\Lambda_{4k+3}-1)(\frac{\Lambda_{4k+3}}{2}-1)\\
&=(2\Lambda_{4(4k+1)}-\Lambda_{4}-\Lambda_{2(4k+1)}+1)(2k\Lambda_{4k+3}+1).\\
\end{align*}
It follows then that $b_{2}(L_{4(4k+1)(4k+3})=2k+1$.

Therefore, $L_{4(4k+1)(4k+3}$ is diffeomorphic to $\#(2k+1) S^2\times S^3$ by Theorem 3.6 and by Corollary 3.9 these manifolds also admit Lorentzian Sasaki-Einstein metrics for all $k\geq 1.$

\noindent\textbf{Case 3}
For the cases $b_2=4,6$ the hypersurfaces of the form
$$\mathcal{Z}=(x^{a_{0}}y+y^{a_{1}}z+z^{a_{2}}t+t^{a_{3}}x=0) $$
will be used. These hypersurfaces were explored by Koll\'ar in \cite{Kol4} which led to examples of rational
homology $\mathbb{CP}^{2}$'s. This class of hypersurfaces is quasi-smooth and well-formed by Theorem 39 of \cite{Kol4}.
Therefore, there are no branch divisors present and therefore the corresponding links over
the well-formed hypersurfaces have no torsion in $H_{2}(L_f,\mathbb{Z})$ by Theorem 3.1.\\
\indent Consider the hypersurface $$\mathcal{Z}_{1213}=(x^{4}y+y^{7}z+z^{10}t+t^{13}x=0)\subset\mathbb{P}(264,157,114,73)$$
of degree $d=1213$. In this case the hypersurface is well-formed by Definition 2.5. Therefore, there are no
branch divisors in $\mathcal{Z}_{1213}$ and so on
the corresponding link $L_{1213}$, $H_{2}(L_{1213},\mathbb{Z})$ is torsion free.
Furthermore, the Betti number is indeed what we want. Observe,
\begin{align*}
divisor\Delta(t)&=\left(\frac{\Lambda_{1213}}{264}-1\right)\left(\frac{\Lambda_{1213}}{157}-1\right)
\left(\frac{\Lambda_{1213}}{73}-1\right)\left(\frac{\Lambda_{1213}}{114}-1\right)=3\Lambda_{1213}+1.\\
\end{align*}
Therefore, $b_{2}(L_{1213})=4$ and by Theorem 3.6 $L_{1213}$ is diffeomorphic to $\#4 S^{2}\times S^{3}$  and it admits
a Sasaki $\eta$-Einstein and Lorentzian Sasaki-Einstein structures since $d-\sum_iw_i>0$. This completes the case $b_2=4.$

Finally, to get $b_2=6$, take the following weighted hypersurface $$\mathcal{Z}_{4435}=(x^{6}y+y^{11}z+z^{16}t+t^{21}x=0)\subset\mathbb{P}(676,379,266,179)$$
of degree $d=4435$. In this case too, $\mathcal{Z}_{4435}$ is well-formed and thus there are no branch divisors which then implies on the corresponding link $L_{4435}$,
$torsH_{2}(L_{4435},\mathbb{Z})=0$. Furthermore,
\begin{align*}
divisor\Delta(t)&=\left(\frac{\Lambda_{4435}}{676}-1\right)\left(\frac{\Lambda_{4435}}{379}-1\right)\left(\frac{\Lambda_{4435}}{266}-1\right)
\left(\frac{\Lambda_{4435}}{179}-1\right)=5\Lambda_{4435}+1
\end{align*}
It follows that $b_{2}(L_{4435})=6$. Since $d-\sum_{i}w_i=2935>0$ we get that $L_{4435}$ admits Sasaki $\eta$-Einstein and Lorentzian
Sasaki-Einstein metrics on $\#6 S^{2}\times S^{3}.$ The proof is complete.$\Box$\\
\indent We have some corollaries. Combining our main theorem together with the work of Boyer, Galicki, Koll\'ar and Nakamaye (\cite{BG4},\cite{BGN1},\cite{BGN3}), we have
\begin{corollary}
The manifolds $\#k S^{2} \times S^{3}$ admit both Sasaki-Einstein
and Lorentzian Sasaki-Einstein structures for all $k\geq 0$.
\end{corollary}
Recall in \cite{Leb}, Catanese and LeBrun were able to show that there exists smooth compact manifolds of dimension $4k$, $k\geq2$, that admit
K\"ahler-Einstein metrics whose scalar curvatures have opposite signs. Since Sasaki $\eta$-Einstein manifolds with constant scalar curvature $2n(\lambda +1)$ are the odd-dimensional analogs of K\"ahler-Einstein metrics and are organized according to $\lambda =-2,>-2,<-2$
one can search for an analogous result of Catanese and LeBrun. Indeed, for our case $n=2$ we get
\begin{corollary}
The manifolds $\#k S^{2} \times S^{3}$ admit pairs of Sasaki $\eta$-Einstein structures $(\mathcal{S}_{1},g_1)$ and $(\mathcal{S}_{2},g_2)$
such that $g_1$ has constant scalar curvature $>-4$ and $g_2$ has constant negative scalar curvature $<-4$ for all $k$.
\end{corollary}
In fact, a $5$-manifold that admits all three types of Sasaki $\eta$-Einstein structures is completely determined.
\begin{corollary} $M^5$ admits a positive, negative and null Sasaki $\eta$-Einstein structure iff $M^{5}$ is diffeomorphic to $\#k S^{2} \times S^{3}$
for $k=3,...,21$.
\end{corollary}
\begin{proof}
Null Sasaki $\eta$-Einstein structures on 5-manifolds are diffeomorphic to k-fold connected sums
of $S^{2}\times S^{3}$ by (\cite{BGM}, \cite{Cuadros}, \cite{Kol1}) for $k=3,...,21.$
\end{proof}
Our main theorem also gives rise to manifolds admitting so-called twistor spinors \linebreak(cf.\cite{Alekseevsky},\cite{Bau},\cite{Bohle}). Recall these are solutions to the conformally invariant
twistor equation on spinors. In fact,  Proposition 3.2 of \cite{Bau} asserts that a simply-connected Lorentzian Sasaki-Einstein manifold admits a twistor spinor and so we arrive at
\begin{corollary}
The manifolds $\#k S^{2} \times S^{3}$ admit a twistor spinor
for $k\geq 0$.
\end{corollary}

\section*{acknowledgments}
I would like to thank Charles Boyer for the encouragement and suggestions in preparing this article. I also would like to thank
the Efroymson Fund for partial financial support under which some of this article was prepared.

\end{document}